\input amstex

\magnification=1150

\centerline{\bf A Combinatorial Curvature Flow for Compact 3-Manifolds with Boundary}

\medskip
\centerline{\bf Feng Luo}
\medskip
\centerline{\bf abstract}
\medskip
We introduce a combinatorial curvature flow for PL metrics on compact triangulated 3-manifolds with boundary consisting of surfaces
of negative Euler characteristic. The flow tends to find the complete hyperbolic metric with totally geodesic boundary on a
manifold. Some of the basic properties of the combinatorial flow are established. 
The most important ones is that the evolution of the combinatorial curvature satisfies a combinatorial heat 
equation. It implies that the total curvature decreases along the flow. The local convergence of the flow to 
the hyperbolic metric is also established if the triangulation is isotopic to a totally
geodesic triangulation.


\medskip
\noindent
\S1. {\bf Introduction}

\medskip
\noindent
1.1.  The purpose of this paper is to announce  the construction of  a combinatorial curvature flow which is a 3-dimensional 
analogy of the work of [CL]. In [CL], we introduced a 2-dimensional combinatorial curvature flow for triangulated
surfaces of non-positive Euler characteristic.  It is shown that for any initial choice
 of PL metric of circle packing type, the flow exists for all time and converges exponentially fast to the
Andreev-Koeb-Thurston's circle packing metrics. 

The basic building blocks for the 2-dimensional flow  are hyperbolic and Euclidean triangles.
In the 3-dimensional case, the basic building blocks are the \it hyperideal tetrahedra \rm discovered by
Bao and Bonahon [BB]. Given an ideal triangulation of a compact 3-manifold with boundary consisting of
surfaces of negative Euler characteristic, we replace each (truncated) tetrahedron by a hyperideal tetrahedron by
assigning the edge lengths. The isometric gluing of these tetrahedra gives a hyperbolic cone metric on the 3-manifold.
The \it PL curvature \rm of the cone metric at an edge is  $2\pi$ less than the sum of dihedral angles at the edge.
The combinatorial curvature flow that we propose is the following system of ordinary differential equations,
$$ d x_i/dt = K_i \tag 1.1$$
where $x_i$ is the length of the i-th edge and $K_i$ is the curvature of the cone metric $(x_1, ..., x_n)$ at the i-th edge.
The equation (1.1) captures the essential features of the 2-dimensional combinatorial Ricci flow in [CL].
The most important of all is that the PL curvature evolves according to a combinatorial heat equation. Thus the corresponding
maximal principle applies. The flow has
the tendency of finding the complete hyperbolic metric of totally geodesic boundary on the manifold.
By analyzing the singularity formations in equation (1.1), it is conceivable that  one could give a new proof of
Thurston's geometrization theorem for these manifolds using (1.1). Furthermore, the flow (1.1) will be a useful tool to find algorithmically 
the complete hyperbolic metric.

\medskip
\noindent
1.2. We now provide some details of the approach.  Suppose $M$ is a compact 3-manifold with 
boundary consisting of surfaces of negative Euler characteristic. An \it ideal triangulation \rm  (or truncated triangulation)
of $M$ is
 the following finite set of data. Take a finite collection of 3-simplexes
and identify their faces in pairs by homeomorphisms. 
The quotient space with a small regular neighborhood of each vertex removed
is homeomorphic to the 3-manifold $M$. By abuse the use of language, we will call (the homeomorphism images of) 
i-dimensional cells ($i \geq 1$) of
the ideal triangulation (truncated triangulation) in the interior of $M$, the \it edges, triangles and tetrahedra. \rm 
A \it hyperideal tetrahedron \rm in the 3-dimensional hyperbolic space is a compact convex polyhedron so that
it is diffeomorphic to a truncated tetrahedron in the 3-dimensional Euclidean space and its four hexagonal  
faces are right-angled hyperbolic hexagons. (Note that
two compact subsets of $\bold R^n$ are diffeomorphic if there is a diffeomorphism between two open neighborhoods of
them sending one compact set to the other). See the figure below.   In the beautiful work of [BB], Bao and Bonahon give a complete characterization
of hyperideal convex polyhedra. As a consequence of their work, hyperideal
tetrahedra are completely characterized  by their six dihedral angles at the six edges so that the sum of three dihedral angles associated
to edges adjacent to each vertex is less than $\pi$.

\medskip

\input epsf
\centerline{\epsfbox{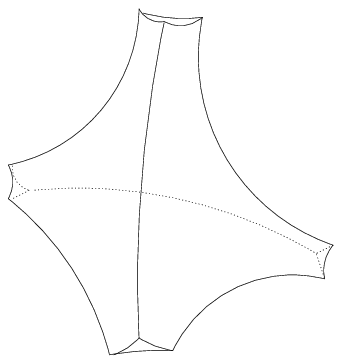}}
\centerline{\it The six edges of a hyperideal tetrahedron are the intersections of its hexagonal faces \rm}
\centerline{\bf figure }

\medskip
One of our main technical observations is the following,

\medskip
\noindent
{\bf Theorem 1.} \it The volume of a hyperideal tetrahedron is a strictly 
concave    function of its dihedral angles. \rm
\medskip

Note that each hyperideal tetrahedron is also determined by its six edge lengths. By the Schlaefli formula, an equivalent statement to
theorem 1 is that the Hessian matrix of the map sending six edge lengths to the six dihedral angle is strictly positive definite.
This  positive definite matrix provides a basis for constructing the combinatorial Laplacian operator for the curvature evolution equation.

\medskip
Now given an ideal triangulated 3-manifold $(M, T)$, let  $E$ be the set of edges in the triangulation and let $n$ be the number
of edges in $E$.  An assignment $x: E \to \bold R_{>0}$ is called a \it hyperbolic cone metric associated to the triangulation
$T$ \rm if for each tetrahedron $t$ in $T$ with edges $e_1, ..., e_6$, the six numbers $x_i = x(e_i)$ (i=1,..., 6) are the
edge lengths of a hyperideal tetrahedron in $\bold H^3$. The set of all hyperbolic cone metrics associated to $T$ is
denoted by $L(M, T)$ which will be considered as an open subset of $\bold R^n = \bold R^E$. The PL curvature of a cone metric
$x \in L(M,T)$ is the map $K(x): E \to \bold R$ sending an edge $e$ to the PL curvature of $x$ at the edge $e$. Again we 
identify the
set of all PL curvatures $\{ K(x) | x \in L(M, T)\}$ with a subset of $\bold R^n$. The
combinatorial curvature flow is the vector field in $L(M, T)$ defined by equation (1.1) where $K_i$ in the
right-hand side is the PL curvature of the metric $x=(x_1,...,x_n)$ at time $t$ at the i-th edge.

\medskip
\noindent
{\bf Theorem 2.} \it  For any ideal triangulated 3-manifold, under the combinatorial curvature flow (1.1), the PL curvature
 $K_i(t)$ evolves according to a combinatorial heat 
equation, 
$$ dK_i(t)/dt = \sum_{j=1}^n a_{ij} K_j(t) \tag 1.2$$
where the matrix $[a_{ij}]_{n \times n}$ is  symmetric negative definite.
 In particular, the total curvature
$\sum_{i=1}^n K_i^2(t)$ is strictly decreasing along the flow unless $K_i(t) =0$ for all $i$.
\rm

\medskip
\noindent
{\bf Theorem 3.} \it For any ideal triangulated 3-manifold $(M,T)$, 
the equilibrium  points  of the combinatorial curvature flow (1.1) are the complete hyperbolic
metric with totally geodesic boundary. Furthermore, each  equilibrium  point is a local attractor of the flow.\rm

\medskip

Another consequence of the convexity is the following local rigidity result for hyperbolic cone metrics without constrains on
 cone angles. Note that  by [HK],  hyperbolic cone metrics with cone angles at most $2\pi$ are  locally rigid.

\medskip
\noindent
{\bf Theorem 4. } \it For any ideal triangulated 3-manifold $(M,T)$, the curvature map $\Pi: L(M, T) \to \bold R^n$ sending
a metric $x$ to its curvature $K(x)$ is a local diffeomorphism. In particular, a hyperbolic cone metric associated to
an ideal triangulation is locally determined by its cone angles.  \rm

\medskip
\noindent

In the rest of the note, we will sketch briefly the ideas of the proofs of the results. In the last section, we propose several 
questions related to the combinatorial curvature flow whose resolution will lead to a new proof of Thurston's geometrization theorem
for this class of 3-manifolds. 

\medskip
\noindent
1.3. Acknowledgment.  We thank Ben Chow for conversations on the topic of Ricci flow and X.S. Lin for 
discussions. The work is supported in part by the NSF and a research grant from the Rutgers University.

\medskip
\noindent
\S2. {\bf Sketch of the Proofs}
\medskip
\noindent
2.1. The key step in the proof of all theorems is to establish the convexity theorem 1. Here is a quick way to see it. 
Suppose $x_1, ..., x_6$ are the lengths of the six edges of a hyperideal tetrahedron so that
the corresponding dihedral angles are $a_1, ..., a_6$. Let $a=(a_1, ..., a_6)$, $x=(x_1, ..., x_6)$ and $V$ be the
volume of the hyperideal tetrahedron. Then
by using cosine laws for hyperbolic triangles and hyperbolic right-angled hexagons, we see that both
functions $x=x(a)$ and $a=a(x)$ are smooth in $a$ and $x$ respectively.  In particular, they are diffeomorphisms.  By the Schlaefli formula $\partial V/\partial a_i = -x_i/2$,
 the Jacobian matrix $[\partial x_i /\partial a_j]_{6 \times 6}$ is
symmetric.  The Jacobian matrix  is non-singular due to the fact that the map $a=a(x)$ is a diffeomorphism. Thus its signature is
independent of the choice of hyperideal tetrahedra.  One checks directly that in the case the hyperideal tetrahedron is regular (i.e., all $x_i$'s are the same and
all $a_i$'s are the same), the Jacobian matrix is positive definite. Thus theorem 1 follows. Since the inverse of
a symmetric positive definite matrix is again symmetric and positive definite, we see that the matrix
 $[\partial a_i/\partial x_j]_{6 \times 6}$ is also symmetric and positive definite.

\medskip
\noindent
2.2. To prove theorem 2, one uses the equation (1.1). Namely,  to understand the evolution of the curvature, it suffices to
understand the evolution of the individual dihedral angle $a_i(t)$. By the chain rule, we have,
$$ da_i/dt = \sum_{j=1}^6 (\partial a_i/\partial x_j) dx_j/dt = \sum_{j=1}^6 (\partial a_i/\partial x_j) K_j. \tag 2.1$$
By the discussion above, the matrix $[\partial a_i/\partial x_j]_{6 \times 6}$ is positive definite.
By summing over all dihedral angles adjacent to a fixed edge, we obtain theorem 2 since the sum of semi-negative definite matrices
are still semi-negative definite. A further study shows  that the resulting matrix is strictly negative definite.

\medskip
\noindent
2.3.  Both theorems 3 and 4 follow from the fact that the combinatorial curvature flow (1.1) is the negative gradient flow of a
strictly convex function defined on the space $L(M, T)$. Indeed, by the Schlaefli formula $dV = -1/2 \sum_i x_i da_i$, 
if we form $H = 2V +\sum_{i=1}^6 a_i x_i$, then
$\partial H/\partial x_i = a_i$. Thus  the function $H$ is a strictly convex function of the edge lengths $x_1, ..., x_6$.
For any hyperbolic cone metric $x$ associated to the ideal triangulation $T$, we define
$$H(x) = 2 vol(x) - \sum_{i} K_i x_i \tag 2.2$$
where $vol(x)$ is the volume of the metric. One checks easily that $H(x)$ is a strictly convex function of $x$ and
the gradient of $H(x)$ is $-(K_1, ..., K_n)$. Thus theorem 3 follows from the standard theory of  ordinary differential equation. To see theorem 4, one uses the fact that the map sending a point to the gradient of a strictly smooth convex function is a local diffeomorphism.

\medskip
\noindent
\S 3. {\bf Some Remarks and Questions}
\medskip
This work and [CL] are motivated by the work of Richard Hamilton [Hi] on Ricci flow. The  strategy  of 
[Hi] is to find a flow deforming metrics so that its curvature evolves according to a heat-type equation. The main focus of
study then shifts from the evolution of the metrics to the evolution of its curvature using the maximal principle for heat equations.
As long as one has control of the curvature evolution, then one gets some control of the metric evolutions by study either the
singularity formation or the long time convergence. Theorem 2 above seems to indicate that the combinatorial flow (1.1) deforms
the cone metrics in the "right" direction. There remains the task of understanding the singularity formations in (1.1) which corresponds
to the degeneration of the hyperideal tetrahedra. This is being investigated. Below are some thoughts on this topic. 
The motivations come from the work of [Th],  [CV], [Riv2], [Lei] and [CL].
 
\medskip
\noindent
3.1.  To understand the singularity formation, we will focus our attention to the function $H$ in (2.2). 
Following [CV] and [Riv2], our goal is to find the (linear) conditions on the ideal triangulation so that it will guarantee the
existence of critical points for $H$. The existence of the ideal triangulation satisfying the (linear) condition will be
related to the topology of the 3-manifold and will be resolved by topological arguments.

Suppose $(M, T)$ is an ideal triangulated compact 3-manifold so that  each boundary component
of $M$ has negative Euler characteristic. A pair $(e, t)$ where $e$ is an edge and $t$ is a tetrahedron containing $e$ is
called a \it corner \rm in $T$.  Following Rivin [Riv1], Casson and Lackenby [La1], we say the triangulated manifold $(M, T)$ 
supports a \it linear hyperbolic structure \rm if one can assign
each corner of $T$ a positive number called the \it dihedral angle \rm so that (1) the sum of  dihedral angles
of all  corners adjacent to each fixed edge
is $2\pi$, and (2) the sum of  dihedral angles of every triple of corners $(e_1, t), (e_2, t), (e_3, t)$ where $e_1, e_2, e_3$
are adjacent to a fixed vertex is less than $\pi$.  By the work of [BB], given a linear hyperbolic structure on $(M, T)$, we can realize
each individual tetrahedron by a hyperideal tetrahedron whose dihedral angles are the assinged numbers so that the sum of
the dihedral angles at each edge is $2\pi$.
It can be shown  that if $(M, T)$ supports a linear hyperbolic structure,
then the manifold $M$ is irreducible without incompressible tori. One would ask if the converse is also true.  The work of Lackenby [La2] gives
some evidences that the following may have a positive answer. See also the work of [KR].

\medskip
\noindent
{\bf Question 1.} \it Suppose $M$ is a compact irreducible 3-manifold with incompressible boundary consisting of surfaces
of negative Euler characteristic. If $M$ contains no incompressible tori and annuli,  is there any ideal triangulation of $M$ which supports a
linear hyperbolic structure? \rm

\medskip
The next question is the 3-dimensional analogy of the 2-dimensional singularity formation analysis in the work of [CV],  [Riv2] 
and [Le]. 
\medskip
\noindent
{\bf Question 2.} \it Suppose $(M, T)$ is an ideal triangulated 3-manifold which supports a  linear hyperbolic structure.
 Does  $H$ have a local minimal point in the space $L(M,T)$ of all cone metrics associated to $(M, T)$?
\rm

\medskip
Positive resolutions of above two questions  will produce a new proof of Thurston's geometrization theorem for this
class of 3-manifolds.


\medskip
\noindent
3.2. Suppose $(M,T)$ supports a linear hyperbolic structure. We define the
volume  of a linear hyperbolic structure to be the sum of the volumes of
its hyperideal tetrahedra. 
Let
$LH(M,T)$ be the space of all linear hyperbolic structures on $(M,T)$.
It can be shown, using Lagrangian multipliers, that the volume function is
strictly concave on $LH(M,T)$ whose maximal point is exactly the complete
hyperbolic metric on $M$. The situation is the same as ideal triangulations
of compact 3-manifolds whose boundary consists of tori. In this case, one realizes each tetrahedron by an ideal tetrahedron in the hyperbolic space. It can 
be shown that
the complete hyperbolic metric of finite volume is exactly equal to the
maximal point of the volume function defined on the space of all linear hyperbolic structures defined in [Ri1]. This was also observed by Rivin [Ri3].
\medskip
\noindent
3.3. We remark that the moduli space of all hyperideal tetrahedra parametrized by their edge lengths $x_1, ..., x_6$ is not
a convex subset of $\bold R^6$. This is the main reason that we have only local convergence and local rigidity in theorems
3 and 4.
However, it is conceivable that
theorem 4 may still be true globally. We do not know yet if the
space  $L(M, T)$ of all cone metrics associated to the ideal triangulated manifold is homeomorphic to a Euclidean space.
It is likely to be the case. In fact, one would hope that there is a diffeomorphism $h: \bold R_{>0} \to \bold R_{>0}$ so
that if we parameterize the space of all hyperideal tetrahedra by $(t_1, ..., t_6)=(h(x_1), ..., h(x_6))$, then 
the space becomes convex in $t$-coordinate. Evidently, if this holds, it implies the space $L(M, T)$ is convex
in the t-coordinate.

\medskip
\medskip

\centerline{\bf References}
\medskip

\noindent
[BB] Bao, Xiliang; Bonahon, Francis: Hyperideal polyhedra in hyperbolic 3-space. Bull. Soc. Math. France 130 (2002), no. 3, 457--491.

\noindent
[CL] Chow, Bennett; Luo, Feng: Combinatorial Ricci flows on surfaces. J. Differential Geom. 63 (2003), no. 1, 97--129.

\noindent
[CV] Colin de Verdière, Yves: Un principe variationnel pour les empilements de cercles.  Invent. Math. 104 (1991), no. 3, 655--669.

\noindent
[Hi] Hamilton, Richard S. Three-manifolds with positive Ricci curvature. J. Differential Geom. 17 (1982), no. 2, 255--306. 

\noindent
[HK] Hodgson, Craig D.; Kerckhoff, Steven P:  Rigidity of hyperbolic cone-manifolds and hyperbolic Dehn surgery. J. Differential Geom. 48 (1998), no. 1, 1--59.

\noindent
[KR] Kang, Ensil and JH Rubinstein, J. H:  ideal triangulations of 3-manifolds I, preprint.

\noindent
[La1] Lackenby, Marc: Word hyperbolic Dehn surgery. Invent. Math. 140 (2000), no. 2, 243--282. 

\noindent
[La2] Lackenby, Marc: Taut ideal triangulations of 3-manifolds. Geom. Topol. 4 (2000), 369--395 (electronic). 

\noindent
[Le]  Leibon, Gregory: Characterizing the Delaunay decompositions of compact hyperbolic surfaces. Geom. Topol. 6 (2002), 361--391 (electronic).

\noindent
[Ri1] Rivin, Igor: Combinatorial optimization in geometry. Adv. in Appl. Math. 31 (2003), no. 1, 242--271.

\noindent
[Ri2] Rivin, Igor:  Euclidean structures on simplicial surfaces and hyperbolic volume. Ann. of Math. (2) 139 (1994), no. 3, 553--580. 

\noindent
[Ri3] Rivin, Igor: private communication.

\noindent
[Th] Thurston, William: Topology and geometry of 3-manifolds, Lecture notes, Princeton University, 1978.
\medskip

Department of Mathematics, Rutgers University, Piscataway, NJ 07059

Email: fluo\@math.rutgers.edu

\end